\documentclass[graybox]{svmult}

\usepackage{mathptmx}       
\usepackage{helvet}         
\usepackage{courier}        
\usepackage{type1cm}        
\usepackage{makeidx}         
\usepackage{graphicx}        
\usepackage{multicol}        
\usepackage[bottom]{footmisc}

\usepackage{amsmath}        
\usepackage{amsfonts}       
\usepackage{amssymb}        
\usepackage[ruled]{algorithm}
\usepackage{algpseudocode}

\makeindex             

\begin{document}

\title*{Parallel implementation of Multilevel BDDC}
\author{Jakub~\v{S}\'{\i}stek, Jan~Mandel, Bed\v{r}ich~Soused\'{\i}k, Pavel Burda}
\authorrunning{Jakub~\v{S}\'{\i}stek et al.}
\institute{
Jakub \v{S}\'{\i}stek 
\at 
Institute~of~Mathematics, Academy~of~Sciences of the~Czech~Republic, \\
\v{Z}itn\'{a}~25, Praha 1, CZ-115 67, Czech Republic, 
\email{sistek@math.cas.cz}
\and 
Jan Mandel
\at
Department of Mathematical and Statistical Sciences, University of Colorado Denver, \\
Campus Box 170, Denver, CO 80217-3364, United States,
\email{jan.mandel@ucdenver.edu}
\and 
Bed\v{r}ich Soused\'{\i}k
\at
Department of Aerospace and Mechanical Engineering, University of Southern California, \\
Olin Hall 430, Los Angeles, CA 90089-2531, United States,
\email{sousedik@usc.edu}
\and 
Pavel Burda
\at
Department of Mathematics, Faculty of Mechanical Engineering, Czech Technical University, \\
Karlovo n\'{a}m\v{e}st\'{\i} 13, Praha 2, CZ-121 35, Czech Republic,
\email{pavel.burda@fs.cvut.cz}
}

%
%
\maketitle

\abstract*{
In application of the Balancing Domain Decomposition by Constraints (BDDC) to a case with many substructures, 
solving the coarse problem exactly becomes the bottleneck which spoils scalability of the solver.
However, it is straightforward for BDDC to substitute the exact solution of the coarse problem by another step of BDDC method
with subdomains playing the role of elements. 
In this way, the algorithm of three-level BDDC method is obtained.
If this approach is applied recursively, multilevel BDDC method is derived.
We present a detailed description of a recently developed parallel implementation of this algorithm.
The implementation is applied to an engineering problem of linear elasticity and a benchmark problem of Stokes flow in a cavity.
Results by the multilevel approach are compared to those by the standard (two-level) BDDC method.}

\abstract{
In application of the Balancing Domain Decomposition by Constraints (BDDC) to a case with many substructures, 
solving the coarse problem exactly becomes the bottleneck which spoils scalability of the solver.
However, it is straightforward for BDDC to substitute the exact solution of the coarse problem by another step of BDDC method
with subdomains playing the role of elements. 
In this way, the algorithm of three-level BDDC method is obtained.
If this approach is applied recursively, multilevel BDDC method is derived.
We present a detailed description of a recently developed parallel implementation of this algorithm.
The implementation is applied to an engineering problem of linear elasticity and a benchmark problem of Stokes flow in a cavity.
Results by the multilevel approach are compared to those by the standard (two-level) BDDC method.}

\section{Introduction}

The Balancing Domain Decomposition by Constraints\index{BDDC} (BDDC) method introduced in \cite{Dohrmann-2003-PSC} 
is one of the most advanced methods of iterative substructuring \index{iterative substructuring} for the solution of large systems of linear algebraic
equations arising from discretization of boundary value problems. 
However, in the case of many substructures, solving the coarse problem exactly becomes 
the limiting factor for scalability of the otherwise perfectly parallel algorithm. 
This has been observed also for the FETI-DP method (e.g. in \cite{Klawonn-2010-HSP}), 
which is closely related to BDDC.
For this reason, recent research in the area is directed towards inexact solutions of the coarse problem.
For example, algebraic multigrid is used in \cite{Klawonn-2010-HSP} to obtain an approximate coarse correction within the FETI-DP method,
and excellent scalability is achieved.

We follow a different approach in this contribution. 
As was mentioned already in \cite{Dohrmann-2003-PSC}, 
it is quite straightforward for BDDC to substitute the exact solution of the coarse problem by another step of the BDDC method
with subdomains playing the role of elements. 
In this way, the algorithm of three-level BDDC method is obtained (studied in \cite{Tu-2007-TBT3D}).
If this step is repeated recursively, one arrives at the \emph{multilevel BDDC method}\index{multilevel BDDC} 
(introduced in \cite{Mandel-2008-MMB} without a parallel implementation).
Unlike for most other domain decomposition methods, 
such extension is natural for BDDC, since the coarse problem has the same structure as the original problem. 
Although the mathematical theory in \cite{Mandel-2008-MMB} suggests worsening of the efficiency of the multilevel BDDC preconditioner\index{preconditioner} with each
additional level,
the resulting algorithm may outperform the standard method with respect to computational time due to better scalability. 
This fact makes the algorithm a good candidate for using on future massively parallel systems.

In this paper, we present a recently developed parallel implementation of multilevel BDDC method.
It is applied to an engineering problem of linear elasticity\index{linear elasticity} and a benchmark problem of Stokes flow\index{Stokes flow} in a cavity.
The results suggest which drawbacks of the two-level implementation might be overcome by the extension to more levels.
Our solver library has been released as an open-source package.

\section{BDDC preconditioner with two and more levels}

The starting point for BDDC is
the \emph{reduced interface problem}  
\ $\mathbf{\widehat{S}} \, \mathbf{\widehat{u}} = \mathbf{\widehat{g}}$,
where $\mathbf{\widehat{S}}$ is the \emph{Schur complement} with respect to \emph{interface}, i.e. unknowns shared by more than one subdomain,
$\mathbf{\widehat{u}}$ is the part of vector of coefficients of finite element basis functions at the interface, 
and $\mathbf{\widehat{g}}$ is sometimes called \emph{condensed right hand side}.
This problem is solved by a Krylov subspace method in the framework of \emph{iterative substructuring}.
Within these methods, application of $\mathbf{\widehat{S}}$ to a vector is realized by parallel solution of 
independent \emph{discrete Dirichlet problems}.
In this way, the costly explicit construction of the Schur complement is avoided.
However, since it is not the main concern of this contribution, the reader is referred to paper \cite{Sistek-2011-APB}, 
or monograph \cite{Toselli-2005-DDM} for details of iterative substructuring.

In what follows, we turn our attention towards the second key part of Krylov subspace methods -- the \emph{preconditioner},
which is realized by one step of the BDDC method.
Let us begin with description of the standard (two-level) version of BDDC.
Let $\mathbf{K}_i$ be the local subdomain matrix, 
obtained by the sub-assembling of element matrices of elements contained in $i$-th subdomain.
We introduce the \emph{coarse space basis functions} on each subdomain represented by columns of matrix ${\varPsi} _i$,
which is the solution to the saddle point problem with multiple right hand sides
\begin{equation}
\label{eq:psi_def}
\left[
\begin{array}
[c]{cc}
    \mathbf{K}_i & \mathbf{C}_i^T \\
    \mathbf{C}_i & \mathbf{0}
\end{array}
\right]  \left[
\begin{array}
[c]{c}
{\varPsi} _i\\
{\varLambda} _i
\end{array}
\right]  =\left[
\begin{array}
[c]{c}
    \mathbf{0}\\
    \mathbf{I}
\end{array}
\right].
\end{equation}
Matrix $\mathbf{C}_i$ represents constraints on functions ${\varPsi} _i$, 
one row per each.
These constraints enforce continuity of approximate solution at \emph{corners} and/or continuity of more general quantities, 
such as averages over shared subsets of interface (\emph{edges} or \emph{faces}) between adjacent subdomains.
The \emph{local coarse matrix} $ \mathbf{K}_{Ci} = {\varPsi} _i^T \mathbf{K}_i {\varPsi} _i = -{\varLambda} _i$ 
is constructed for each subdomain.
The \emph{global coarse matrix} $\mathbf{K}_{C}$ is obtained by the assembly procedure from local coarse matrices.
This can be formally written as $\mathbf{K}_{C} = \sum _{i = 1}^N \mathbf{R}_{Ci}^T \mathbf{K}_{Ci} \mathbf{R}_{Ci}$,
where $\mathbf{R}_{Ci}$ realize the restriction of global coarse degrees of freedom to local coarse degrees of freedom of $i$-th subdomain.

Suppose 
$\mathbf{\widehat{r}} =\mathbf{\widehat{g}} - \mathbf{\widehat{S}} \, \mathbf{\widehat{u}}$
is a~residual within the Krylov subspace method.
The residual assigned to $i$-th subdomain is computed as $\mathbf{r}_i = \mathbf{E}_i^T \,\mathbf{\widehat{r}}$,
where matrices of weights $\mathbf{E}_i^T$ distribute $\mathbf{\widehat{r}}$ to subdomains.
The subdomain correction is now defined as the solution to the system 
\begin{equation}
\label{eq:sub_cor}
\left[
\begin{array}
[c]{cc}
    \mathbf{K}_i & \mathbf{C}_i^T \\
    \mathbf{C}_i & \mathbf{0}
\end{array}
\right]  \left[
\begin{array}
[c]{c}
\mathbf{z} _i\\
\mathbf{\lambda} _i
\end{array}
\right]  =\left[
\begin{array}
[c]{c}
    \mathbf{r}_i\\
    \mathbf{0}
\end{array}
\right].
\end{equation}
The residual for the coarse problem is constructed using the coarse basis functions subdomain by subdomain and assembling the contributions as 
$ \mathbf{r}_C = \sum _{i = 1}^N \mathbf{R}_{Ci}^T {\varPsi} _i^T \mathbf{E}_i^T \mathbf{\widehat{r}}$.
The coarse correction is defined as the solution to problem
$ \mathbf{K}_C \, \mathbf{z}_C = \mathbf{r}_C$.
Both corrections are finally added together and averaged on the interface by matrices $\mathbf{E}_i$ to produce the preconditioned residual
$\mathbf{\widehat{z}} = \sum _{i = 1}^N \mathbf{E}_i \left( {\varPsi} _i \mathbf{R}_{Ci} \mathbf{z}_C  + \mathbf{z} _i \right)$.

In the three-level BDDC method \cite{Tu-2007-TBT3D},
the matrix $\mathbf{K}_C $ is not constructed on the second level.
Instead, subdomains from the basic (first) level are grouped into subdomains on the next (second) level in the same way as elements of the original mesh are grouped into subdomains of the first level.
The whole procedure described in this section is now repeated for the second level and thus the final coarse problem 
represents the third level.
Obviously, this can be repeated again in the multilevel BDDC method.
The only important difference between the first and the higher levels is the additional \emph{interior pre-correction} and \emph{post-correction} applied on higher levels in order to approximate the whole vector of coarse solution on the lower level.

According to \cite{Mandel-2008-MMB}, the condition number of the operator preconditioned by multilevel BDDC with $L$ levels satisfies
$
\kappa(\mathbf{{M}}_{BDDC}\mathbf{\widehat{S}}) \leq
{\prod_{\ell=1}^{L-1}}
C_{\ell}\left(  1+\log\frac{H_{\ell}}{H_{\ell-1}}\right)  ^{2},
$
where $H_{\ell}$ is the characteristic size of subdomain on level $\ell$,
and $H_{0} \equiv h$ is the characteristic size of element. 
Index $\ell$ is used here and throughout the next section to denote particular level.
Due to the product present in this bound, each additional level worsens the mathematical efficiency of the multilevel preconditioner.
The proof of the condition number bound as well as details of the algorithm of multilevel BDDC can be found in \cite{Mandel-2008-MMB}.

\section{Parallel implementation}

Our implementation of the multilevel BDDC method has been recently released as an open-source solver 
library BDDCML\footnote{\url{http://www.math.cas.cz/$\sim$sistek/software/bddcml.html}}.
It is written in Fortran 95 programming language and parallelized by MPI.
The solver relies on the sparse direct solver MUMPS
--- a serial instance is used for each subdomain problem and a parallel instance is called for the final coarse problem.
The solver supports assignment of several subdomains to each processor, 
since it is often useful to create divisions independently of number of available processors.
A~division of the mesh into subdomains on the first level is either provided to the solver by user's application or created internally by ParMETIS.
The METIS package is currently used for this purpose on higher levels.

Similarly to other related preconditioners, we first need to set-up the multilevel BDDC preconditioner, 
which is then applied in each iteration of the Krylov subspace method.
Details of the set-up are given in Algorithm~\ref{alg:setup}, while key operations of each application are summarized in Algorithm~\ref{alg:prec_appl}.
In these descriptions, we provide comments on how the steps are implemented in BDDCML in parentheses.

\begin{algorithm}[htbp]
    \caption{Set-up of BDDC preconditioner with $L$ levels}
    \label{alg:setup}
    \begin{algorithmic}[1]
    \For {level $\ell =  1, \dots, L-1 $ }
        \If{$\ell > 1$}
           \State build \emph{pseudo-mesh}: subdomains $\rightarrow$ `elements'; corners + edges + faces $\rightarrow$ `nodes'
        \EndIf
        \State divide \emph{pseudo-mesh} into subdomains (by METIS for $\ell > 1$, or by ParMETIS for $\ell = 1$)
        \State classify interface into \emph{faces}, \emph{edges}, \emph{vertices}
        \State select \emph{corners} (using face-based algorithm from \cite{Sistek-2011-FSC})
        \State assemble matrices of subdomains $\mathbf{K}_i^{\ell}$ (use MPI to collect them on assigned cores)
        \State prepare interior correction -- factorize interior block of $\mathbf{K}_i^{\ell}$ (serial MUMPS)
        \State factorize the matrices of local saddle point problems (\ref{eq:psi_def}) (serial MUMPS)
        \State find \emph{coarse basis functions} $\varPsi _i^{\ell}$  and \emph{coarse matrices} $\mathbf{K}_{Ci}^{\ell}=-{\varLambda}_i^{\ell}$ from (\ref{eq:psi_def}) (serial MUMPS)
    \EndFor
    \State factorize \emph{global coarse matrix} $\mathbf{K}_C^{L-1} = \sum _{i = 1}^{N_{L-1}} (\mathbf{R}_{Ci}^{L-1})^T \mathbf{K}_{Ci}^{L-1} \mathbf{R}_{Ci}^{L-1}$ (parallel MUMPS)
    \end{algorithmic}
\end{algorithm}

\begin{algorithm}[htbp]
    \caption{Application of BDDC preconditioner with $L$ levels}
    \label{alg:prec_appl}
    \begin{algorithmic}[1]
    \For {level $\ell =  1, \dots, L-1 $ }
        \If{$\ell > 1$}
            \State $\mathbf{\widehat{r}}^{\ell} \gets \mathbf{r}_C^{\ell - 1}$
            \State compute \emph{interior pre-correction} of residual $\mathbf{\widehat{r}}^{\ell}$ (serial MUMPS)
        \EndIf
        \State distribute residual among subdomains $\mathbf{r}_i^{\ell} = (\mathbf{E}_i^{\ell})^T \,\mathbf{\widehat{r}}^{\ell}$
        \State determine subdomain corrections $\mathbf{z}_i^{\ell}$ from (\ref{eq:sub_cor}) (serial MUMPS)
        \State construct coarse residual $ \mathbf{r}_C^{\ell} = \sum _{i = 1}^{N_{\ell}} (\mathbf{R}_{Ci}^{\ell})^T ({\varPsi} _i^{\ell})^T (\mathbf{E}_i^{\ell})^T \mathbf{\widehat{r}}^{\ell}$ (collective MPI)
    \EndFor
    \State solve the coarse problem $\mathbf{K}_C^{L-1} \, \mathbf{z}_C^{L-1} = \mathbf{r}_C^{L-1}$ (parallel MUMPS)
    \For {level $\ell = L-1, \dots, 1$ }
        \If{$\ell < L-1$}
            \State $\mathbf{z}_C^{\ell} \gets \mathbf{\widehat{z}}^{\ell + 1}$
        \EndIf
        \State combine coarse correction and subdomain corrections  $\mathbf{\widehat{z}}^{\ell} = \sum _{i = 1}^{N_{\ell}} \mathbf{E}_i^{\ell} \left( {\varPsi} _i^{\ell} \mathbf{R}_{Ci}^{\ell} \mathbf{z}_C^{\ell}  + \mathbf{z} _i^{\ell} \right)$
        \If{$\ell > 1$}
            \State apply \emph{interior post-correction} to $\mathbf{\widehat{z}}^{\ell}$ (serial MUMPS)
        \EndIf
    \EndFor
    \end{algorithmic}
\end{algorithm}

\section{Numerical results}

\begin{figure}[htbp]
\sidecaption
\includegraphics[width=0.49\textwidth]{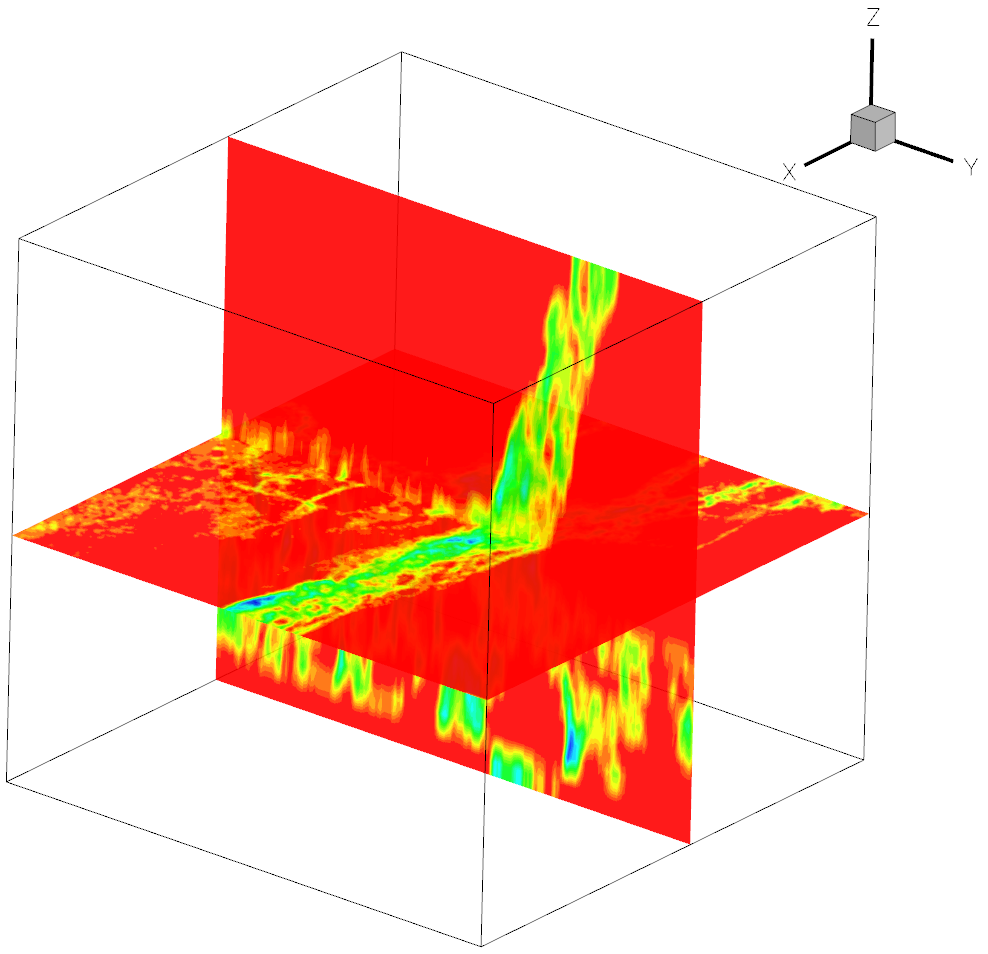}
\includegraphics[width=0.47\textwidth]{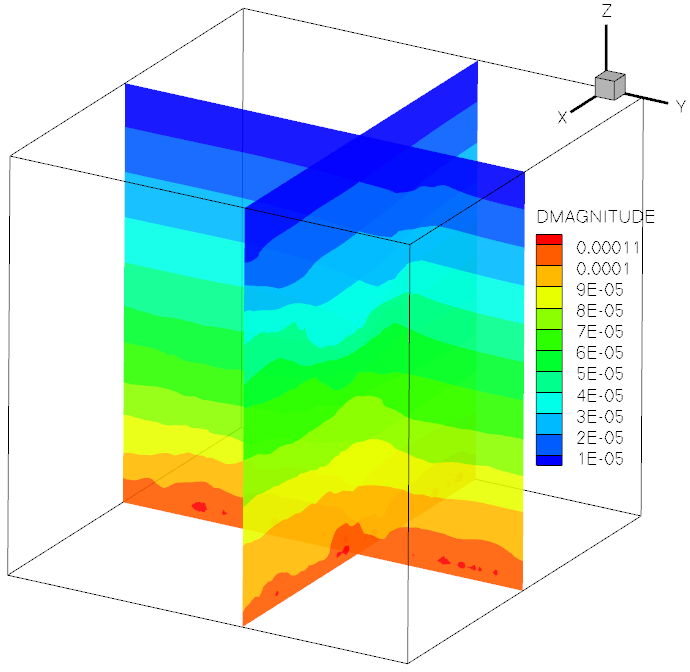}
\caption{Geocomposite problem: slices through material distribution (left) and displacement field (right).}
\label{fig:geoc}
\end{figure}

The first example corresponds to a problem of mechanical analysis of a cubic sample of geocomposite and was analyzed in \cite{Blaheta-2009-SDD}.
The length of the edge of the cube is 75 mm.
The cube comprises five distinct materials identified by means of computer tomography (Fig. \ref{fig:geoc} left),
which causes anisotropic response of the cube even for simple axial stretching in z direction (Fig. \ref{fig:geoc} right).
The problem is discretized using unstructured grid of about 12 million linear tetrahedral elements,
resulting in approximately 6 million unknowns.
The mesh was divided into 1024, 128, and 16 subdomains on the first, second and third level, respectively, 
and the respective coarse problems (using corners and arithmetic averages on all edges and faces) contain 
86,094, 11,265, and 612 unknowns.

Table \ref{tab:geoc_numit} summarizes the efficiency of the multilevel preconditioner 
by means of the resulting condition number (estimated from the tridiagonal matrix generated during iterations of preconditioned conjugate gradient (PCG\index{PCG}) method) 
and number of iterations.
The iterations were stopped when the relative residual $\|\mathbf{\widehat{r}}\|/\|\mathbf{\widehat{g}}\|$ decreased bellow $10^{-6}$.
This table confirms the predicted worsening of the condition number with each additional level expected from the condition number bound.

Table \ref{tab:geoc_scaling} contains a~strong scaling test using different number of levels.
We differentiate the time spent on set-up and in PCG.
All these computations were performed on the IBM SP6 computer at CINECA, Bologna.
The computer is based on IBM Power6 4.7 GHz processors with 4 GB of RAM per core.

We can conclude from Tab. \ref{tab:geoc_scaling} that while adding levels seems not to be feasible for small number of cores (the computational time stagnates or even grows), 
it improves the scaling on many cores.
The minimal overall solution time is achieved for four levels and largest number of cores,
despite the largest number of required iterations.

\begin{table}
\caption{Condition number and number of iterations for different number of levels.}
\label{tab:geoc_numit}
\begin{tabular}
[c]{p{25mm}p{25mm}p{25mm}p{25mm}}
\hline
num. of levels               &  num. of subs.        &  cond. num.    & num. of PCG its. \\
\hline
2                            &  1024/1                    & 50                 & 46 \\
3                            &  1024/128/1                & 79                 & 56 \\
4                            &  1024/128/16/1             & 568                & 131 \\
\hline
\end{tabular}                                        
\end{table}

\begin{table}
\caption{Strong scaling for geocomposite problem using two, three, and four levels.}
\label{tab:geoc_scaling}
\begin{tabular}
[c]{p{35mm}ccccc}
\hline
number of processors         &  64    & 128   & 256   & 512   &  1024 \\\hline
\multicolumn{6}{l}{\textbf{2 levels}} \\
BDDC set-up time (s)           &  61.0  & 37.7  & 25.7  & 23.2  &  39.5 \\
PCG time (s)             &  22.3  & 19.9  & 27.8  & 44.9  &  97.5 \\\hline
\multicolumn{6}{l}{{\bf 3 levels}} \\
BDDC set-up time (s)           &  49.5  & 29.0  & 18.4  & 12.6  &  11.0 \\
PCG time (s)             &  28.5  & 22.6  & 16.7  & 14.7  &  13.2 \\\hline
\multicolumn{6}{l}{{\bf 4 levels}} \\
BDDC set-up time (s)           &  49.4  & 28.6  & 17.8  & 12.3  &  9.1 \\
PCG time (s)             &  60.6  & 33.2  & 21.2  & 15.4  &  11.8 \\\hline
\end{tabular}                                        
\end{table}

Our second example is a problem of Stokes flow in a 3D lid driven cavity. 
We use the set-up suggested in \cite{Wathen-2003-NPO}: zero velocity is prescribed on all faces of the $[0,1]^3$ cube except the face for z = 1, 
where unit velocity vector $\mathbf{u} = [ 1/\sqrt{3},\sqrt{2/3},0] $ is prescribed.
We have used this test case also in the recent paper \cite{Sistek-2011-APB}, 
but we have not presented parallel results there. 

\begin{figure}[htbp]
\sidecaption
\includegraphics[width=0.43\textwidth]{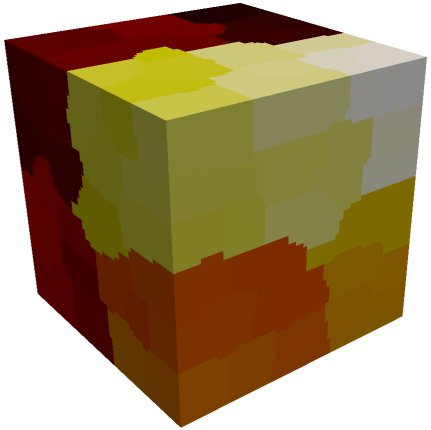}
\hskip 5mm
\includegraphics[width=0.55\textwidth]{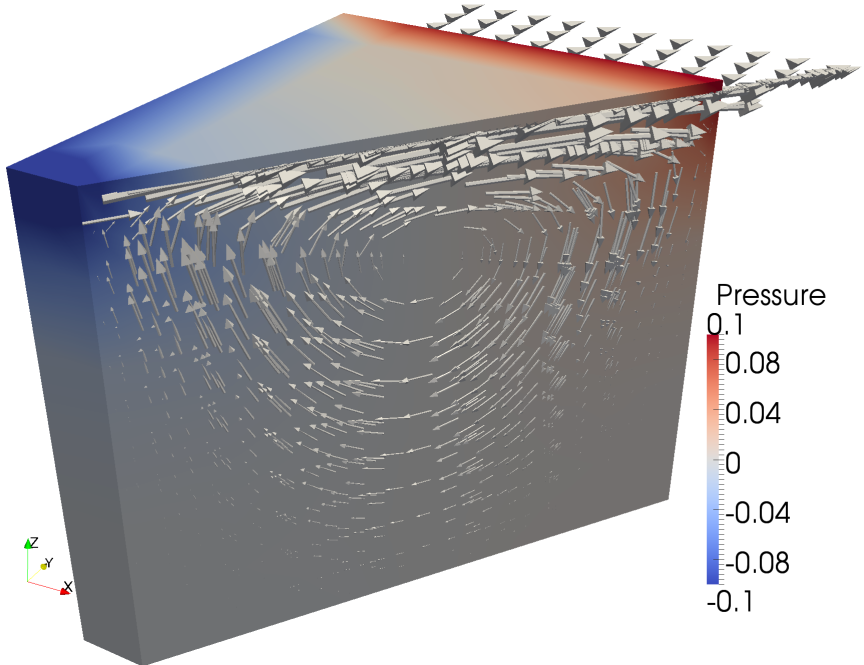}
\caption{Stokes flow in lid driven cavity: example of division into 64 subdomains (left), pressure contours and velocity vectors in the 
cut in direction of the prescribed velocity
$[ 1/\sqrt{3},\sqrt{2/3},0] $ 
(right).}
\label{fig:cavity}
\end{figure}

The problem is uniformly discretized using hexahedral Taylor--Hood finite elements.
Computational mesh was divided into irregular partitions using the METIS graph partitioner (see Fig.~\ref{fig:cavity} for an example).
A~plot of pressure inside the cavity and velocity vectors is given in Fig.~\ref{fig:cavity} (right).
 
Table \ref{tab:scaling_cavity} summarizes a~weak scaling test for this problem. The sequence of problems ranging from 1.7 million unknowns to 25.4 million 
unknowns are distributed among processors such that the size of local problems is kept approximately constant around 50 thousand unknowns.
We present results by the BDDC preconditioner using two and three levels combined with BiCGstab\index{BiCGstab} method.
For comparison, we also report results of the preconditioner based on Additive Schwarz\index{Additive Schwarz} Method (ASM) combined with BiCGstab 
and using MUMPS offered in the PETSc library version 3.1. 
We report numbers of iterations and required overall computational times.
Where `n/a' is present in the table, the simple ASM method was unable to provide the solution, 
which was often encountered for larger problems.
The computational times were obtained on Darwin supercomputer of the University of Cambridge, 
using Intel Xeon 5100  3.0 GHz processors with 2 GB of RAM per core.

We can see in Tab.~\ref{tab:scaling_cavity} that
number of BiCGstab iterations remains almost constant for the two-level method, while mildly growing when using three levels, 
being again larger for the latter.
The computational time slightly grows with problem size for BDDC, both using two and three levels, but this growth is rather acceptable. 
More importantly, we can also see that for this case the benefit of using an additional level is slightly outweighed by the overhead of the additional iterations,
and so the computational time is not improved by using three levels.
Finally, the comparison with ASM results seem to confirm the computational efficiency of our solver which is comparable with the state-of-the-art PETSc library,
and the rather poor performance of ASM in comparison to BDDC in this setting.

\begin{table}
\caption{Weak scaling for Stokes flow in the cavity: additive Schwarz method (ASM), and BDDC using two and three levels.}
\label{tab:scaling_cavity}
\begin{tabular}
[c]{p{9mm}p{9mm}p{12mm}p{7mm}p{14mm}p{7mm}p{14mm}p{12mm}p{7mm}p{12mm}}
\hline 
 \multicolumn{3}{c}{}         & \multicolumn{2}{l}{ASM} & \multicolumn{2}{l}{BDDC (2 levels)} & \multicolumn{3}{l}{BDDC (3 levels)} \\\hline
\# elms.    &   \# dofs.  &   \# cores & \# its.  &    time (s)            &    \# its. &  time (s)                 & divisions &  \# its. &  time (s)  \\\hline
40$^3 $     &   1.7M  &   32    & 282  &    533             &    18      &  122              & 32/4/1   &  22   &  126        \\
50$^3 $     &   3.2M  &   64    & 396  &    805             &    19      &  132              & 64/8/1   &  25   &  205        \\
64$^3 $     &   6.7M  &  128    & 384  &    536             &    21      &  186              & 128/16/1 &  30   &  194        \\
80$^3 $     &  13.1M  &  256    & n/a  &    n/a             &    21      &  178              & 256/32/1 &  36   &  201        \\
100$^3$     &  25.4M  &  512    & n/a  &    n/a             &    20      &  205              & 512/64/1 &  35   &  211        \\
\hline
\end{tabular}                   
\end{table}

\section{Conclusion}

We have presented a~parallel open-source implementation of the multilevel BDDC method.
The two-level algorithm has scalability issues related to the coarse problem solution, mainly in the part of iterations. 
It can be noted that for the tested cases, 
it has not been the size of the coarse problem, but rather its fragmentation among too many cores which causes these issues.
From our experiments, it appears that the multilevel preconditioner tends to scale better in both parts -- set-up and Krylov subspace iterations.
While the better scalability is able to translate into much faster solution for some cases, 
the extra overhead can also just cancel out the savings for other cases. 
It is therefore important to choose appropriate number of levels for a~particular problem.
We expect that advantages of the multilevel approach would pronounce further for problems divided into many (tens of thousands) of subdomains. 
Such challenging problems will likely become common in near future and will provide valuable feedback for further research in this field.

\begin{acknowledgement}
We are grateful to Prof. Blaheta and Dr. Star\'{y} (Institute of Geonics AS CR) for providing the geocomposite problem.
We are also grateful to Dr. Cirak (University of Cambridge) for providing the OpenFTL package 
and computer time on Darwin.
This work was supported by Ministry of Education, Youth and Sports of the Czech Republic under research project LH11004,
by Czech Science Foundation under project 106/08/0403,
by Institutional Research Plan AV0Z 10190503 of the AS CR,
by grant IAA100760702 of the Grant Agency of AS CR,
and by National Science Foundation under grant DMS-0713876. 
The research was started during two visits of Jakub \v{S}\'{\i}stek at the University of Colorado Denver
and some parts of the work have been performed under the HPC-Europa2 project with the support of the European Commission.
\end{acknowledgement}


\end{document}